\documentclass[11pt,reqno, twoside]{article}
\usepackage{amsfonts}
\usepackage{amsfonts, amsmath, amssymb, amscd, amsthm}

\usepackage{mathrsfs}
\usepackage{multicol}
\usepackage{comment}

 \theoremstyle{definition}
 
 \theoremstyle{remark}

 \theoremstyle{claim}
 \numberwithin{equation}{section}

\def\R{\mathbb R}

\def\S{\mathbb S}

\def\f{\frac}

\newcommand{\p}{\partial}

\numberwithin{equation}{section}

\begin{document}
\title{A Gap Theorem for Self-shrinkers of the Mean Curvature Flow
in Arbitrary Codimension\thanks {The first author was partially supported by NSF grant DMS-0909581; the second author was supported by NSFC 10971110.}
}
\author{ Huai-Dong Cao and Haizhong Li}
\date{}
\maketitle

\begin{abstract}
In this paper, we prove a classification theorem for self-shrinkers of the mean curvature flow with $|A|^2\le 1$ in
arbitrary codimension. In particular, this implies a gap theorem for self-shrinkers  in
arbitrary codimension.
\end{abstract}

\section {Introduction}

Let $x:M^n\to \R^{n+p}$ be an $n$-dimensional submanifold in
the (n+p)-dimensional Euclidean space. If we let the position
vector $x$ evolve in the direction of the mean curvature $\bf{H}$,
then it gives rise to a solution to the mean curvature flow:
\begin{equation}\label{1-1}
x:M\times [0,T)\rightarrow \R^{n+p}, \qquad \f{\p x}{\p t} = \bf{H}
\end{equation}

We call the immersed manifold $M$ a self-shrinker if it satisfies
the quasilinear elliptic system:
\begin{equation}\label{1-2}
\bf{H}=-x^{\perp}
\end{equation}
where $\perp$ denotes the projection onto the normal bundle of $M$.

Self-shrinkers are an important class of solutions to the mean curvature flow \eqref{1-1}.
Not only they are shrinking homothetically under mean curvature flow
(see, e.g., \cite{CM}), but also they describe possible blow
ups at a given singularity of the mean curvature flow.

In the curve case, U. Abresch and J. Langer \cite{AL} gave a complete classification of
all solutions to \eqref{1-2}. These curves are so-called Abresch-Langer curves.

In the hypersurface case (i.e. codimension 1), K. Ecker and G. Huisken \cite{EH} proved
that if an entire graph with polynomial volume growth is a self-shrinker, then it  is
necessarily a hyperplane. Recently L. Wang \cite{LW} removed the condition of {\it
polynomial volume growth} in Ecker-Huisken's Theorem. Let $|A|^2$ denote the norm square
of the second fundamental form of $M$. In \cite{H2} and \cite{H3}, G. Huisken proved a
classification theorem that $n$-dimensional self-shrinkers  satisfying (1.2) in
$\R^{n+1}$ with non-nengative mean curvature, bounded $|A|$,  and polynomial volume
growth are $\Gamma\times \R^{n-1}$, or $\S^m (\sqrt{m})\times\R^{n-m}$ ($0\leq m\leq n$).
Here, $\Gamma$ is a Abresch-Langer curve and $\S^m(\sqrt{m})$ is a $m$-dimensional sphere
of radius $\sqrt{m}$. Recently, T.H. Colding and W.P. Minicozzi \cite{CM} showed that G.
Huisken's classification theorem still holds without the assumption that {\it $|A|$ is
bounded.} Moreover, they showed that the only embedded entropy stable self-shrinkers with
polynomial volume growth in $\R^{n+1}$ are hyperplanes, n-spheres, and cylinders.

In arbitrary codimensional case, K. Smoczyk  \cite{Smo} proved the following two results:
(i) For any  $n$-dimensional compact self-shrinker $M^n$ in $R^{n+p}$ satisfying (1.2), if
${\bf H}\not=0$ and unit mean curvature vector field $\nu={\bf H}/|{\bf H}|$ is parallel
in the normal bundle, then $M^n=\S^n(\sqrt{n})$ in $\R^{n+1}$; (ii) For any $n$-dimensional
compact self-shrinker $M^n$ in $R^{n+p}$ satisfying (1.2), if $M^n$ is a complete
self-shrinker with ${\bf H}\not= 0$ and unit mean curvature vector field $\nu={\bf
H}/|{\bf H}|$ is parallel in the normal bundle, and having uniformly bounded geometry,
then $M^n$  is either $\Gamma\times \R^{n-1}$, or $N^m\times \R^{n-m}$. Here $\Gamma$ is
an Abresch-Langer curve and $N^m$ is a $m$-diemnsional minimal submanifold in
$\S^{m+p-1}(\sqrt{m})$. On the other hand, Q. Ding and Z. Wang \cite{DW} recently have extended the
result of L. Wang \cite{LW} to higher codimensional case under the condition of {\it flat
normal bundle}.

Very recently,  based on an identity of Colding and Minicozzi (see (9.42) in \cite{CM}), N. Q. Le and N. Sesum \cite{LS} proved  a gap theorem (cf. Theorem 1.7 in \cite{LS})
for self-shrinkers of codimension 1: if a hypersurface $M^n\subset \R^{n+1}$ is a smooth complete
embedded self-shrinker without boundary and with polynomial volume
growth, and satisfies $|A|^2<1$, then  $M^n$ is a hyperplane. Motivated by this result of Le and Sesum,  we prove in this paper the
following classification theorem for self-shrinkers in arbitrary codimensions:

\medskip\noindent
{\bf Theorem 1.1} {\it If $M^n \to \R^{n+p}$ ($p\ge 1$) is an
$n$-dimensional complete self-shrinker without boundary and with
polynomial volume growth, and satisfies
\begin{equation}
|A|^2\leq 1,
\end{equation}
then $M$ is one of the followings:
\smallskip

(i) a round sphere $\S^n(\sqrt{n})$ in $\R^{n+1}$,

(ii) a cylinder $\S^m(\sqrt{m})\times \R^{n-m}$, $1\leq m\leq n-1$,
in $\R^{n+1}$,

(iii) a hyperplane in $\R^{n+1}$.
\smallskip

\noindent Here $|A|^2$ is the norm square of the second fundamental
form of $M$.}

\medskip

As an immediate consequence, we have the following gap theorem valid for arbitrary codimensions:

\medskip\noindent {{\bf Corollary 1.1} {\it If
$M^n\to \R^{n+p}$ ($p\ge 1$) is a smooth complete
embedded self-shrinker without boundary and with polynomial volume
growth, and satifies
\begin{equation}
|A|^2<1,
\end{equation}
then  $M$ is a hyperplane in $\R^{n+1}$.}

\medskip\noindent {\bf Remark 1.1} We expect that the condition on volume growth in Theorem 1.1 and Corollary 1.1 can be removed.
In fact, it was conjectured by the first author that a complete self-shrinker automatically has polynomial volume
growth. Note that D. Zhou and the first author \cite{CZ} proved that a complete Ricci shrinker necessarily has at most Euclidean
volume growth.

\medskip\noindent {\bf Remark 1.2} Shortly after our work was
finished, Q. Ding and Y. L. Xin \cite{DX} proved that  any
complete non-compact {\it properly immersed} self-shrinker $M^n$ in
$\R^{n+p}$ has at most Euclidean volume growth.

\noindent {\bf Acknowledgements.} Part of the work was carried out while the first author was visiting the Mathematical Sciences Center of Tsinghua University during fall 2010.  He would like to thank the Center for its hospitality and support. The authors would
also like to thank the referee for helpful comments which make the proofs of Lemma 3.1 and Proposition 5.1 more readable.

\section { Preliminaries}

In this section, we recall some formulas and notations for
submanifolds in Euclidean space by using the method of moving frames.

 Let $x:M^n\to \R^{n+p}$ be an  $n$-dimensional  submanifold of the
 $(n+p)$-dimensional Euclidean space $R^{n+p}$. Let $\{e_1,\cdots,e_{n}\}$
be a local orthonormal basis of $M$ with respect to the induced
metric, and $\{\theta_1,\cdots, \theta_n\}$ be their dual 1-forms. Let
$e_{n+1},\cdots, e_{n+p}$ be the local unit orthonormal normal
vector fields.

In this paper we make the following convention on the range of
indices:

$$
1\leq i,j,k\leq n; \qquad
n+1\leq \alpha,\beta,\gamma\leq n+p.
$$

Then we have the following structure equations,
\begin{equation}\label{2-1}
dx=\sum\limits_i \theta_i e_i,
\end{equation}

\begin{equation}\label{2-2}
de_i=\sum\limits_j\theta_{ij}e_j+\sum\limits_{\alpha,j}
h^\alpha_{ij} \theta_je_\alpha,
\end{equation}

\begin{equation}\label{2-3}
de_\alpha=-\sum\limits_{i,j}h^\alpha_{ij}\theta_j e_i+\sum
\limits_\beta \theta_{\alpha\beta}e_\beta,
\end{equation}
where  $h^\alpha_{ij}$ denote the components of the second fundamental
 form of $M$ and $\theta_{ij}$, $\theta_{\alpha\beta}$ denote the connections of the tangent bundle and normal bundle of $M$, respectively.

 The Gauss equations are given by
\begin{equation}\label{al}R_{ijkl}=\sum_{\alpha}
(h_{ik}^{\alpha}h_{jl}^{\alpha}-h_{il}^{\alpha}h_{jk}^{\alpha})
\end{equation}
\begin{equation}\label{am}R_{ik}=\sum_{\alpha}
H^{\alpha}h_{ik}^{\alpha}-\sum_{\alpha,j}h_{ij}^{\alpha}h_{jk}^{\alpha}
\end{equation}
\begin{equation}\label{an}R =H^{2}-|A|^2
\end{equation}
where $R$ is the scalar curvature of $M$,
$|A|^2=\sum\limits_{\alpha,i,j}(h^\alpha_{ij})^2$ is the norm square
of the second fundamental form, ${\bf
H}=\sum\limits_{\alpha}H^\alpha e_\alpha=\sum\limits_{\alpha}
(\sum\limits_i h^\alpha_{ii})e_\alpha$ is the mean curvature vector
field, and $H=|\bf{H}|$ is the mean curvature of $M$.

The Codazzi equations are given by (see, e.g.,  \cite{Li})

\begin{equation}\label{c4}h^\alpha_{ijk}=h^\alpha_{ikj},
\end{equation}
where the covariant derivative of $h^\alpha_{ij}$ is defined by
\begin{equation}\label{2-8}
\sum_{k}h^\alpha_{ijk}\theta_k=dh^\alpha_{ij}+\sum_k
h^\alpha_{kj}\theta_{ki} +\sum_k h^\alpha_{ik}\theta_{kj}+\sum_\beta
h^\beta_{ij}\theta_{\beta\alpha}.
\end{equation}

If we denote by $R_{\alpha\beta i j}$ the  curvature tensor of
the normal connection $\theta_{\alpha\beta}$ in the normal bundle of $x:
M \rightarrow \R^{n+p}$, then the Ricci equations are
\begin{equation}\label{aj}R_{\alpha\beta i j}=\sum_{k}
(h_{ik}^{\alpha}h_{kj}^{\beta}-h_{jk}^{\alpha}h_{ki}^{\beta}).
\end{equation}

By exterior differentiation of $(2.8)$, we have the following Ricci
identities (see, e.g., \cite{Li})
\begin{equation}
h^\alpha_{ijkl}-h^\alpha_{ijlk}=\sum\limits_m h^\alpha_{mj}R_{mikl}
+\sum\limits_m h^\alpha_{im}R_{mjkl}+\sum\limits_\beta h^\beta_{ij}
R_{\beta\alpha kl}.
\end{equation}

We define the first and second covariant derivatives, and Laplacian of the mean curvature vector field ${\bf H}=\sum\limits_\alpha H^\alpha e_\alpha$ in the normal bundle $N(M)$ as follows (cf. \cite{CL}, \cite{Li})

\begin{equation}
\sum\limits_i H^\alpha_{,i}\theta_i=dH^\alpha+\sum\limits_\beta H^\beta\theta_{\beta\alpha},
\end{equation}

\begin{equation}
\sum\limits_j H^\alpha_{,ij}\theta_j=dH^\alpha_{,i}+
\sum\limits_j H^\alpha_{,j}\theta_{ji}+\sum\limits_\beta H^\beta_{,i}\theta_{\beta\alpha},
\end{equation}

\begin{equation}
\Delta^\perp H^\alpha=\sum\limits_i H^\alpha_{,ii},\qquad
H^\alpha=\sum\limits_k h^\alpha_{kk}.
\end{equation}

Let $f$ be a smooth function on $M$, we define the
covariant derivatives $f_i$, $f_{ij}$, and the Laplacian of $f$ as
follows
\begin{equation}\label{2-14}
df=\sum_i f_i\theta_i,\qquad \sum_j
f_{ij}\theta_j=df_i+\sum_jf_j\theta_{ji},\qquad \Delta f= \sum_i
f_{ii}.
\end{equation}

\section {A Key Lemma}
As we mentioned in the introduction, the proof of Le-Sesum's gap theorem relies on an important identity of Colding and Minicozzi \cite{CM}
for hypersurfaces. The identity, see (9.42) in \cite{CM} or (4.1) in \cite{LS}, is obtained in terms of certain second order linear operator for hypersurfaces
(which is part of the Jacobi operator for the second variation).
In this section, we  derive a similar inequality for arbitrary codimensions.

Let $a$ be any fixed vector in $\R^{n+p}$, we define the following
height functions in the $a$ direction on $M$,
\begin{equation}
f=\langle x,a\rangle, 
\end{equation}
and
\begin{equation}
g_\alpha =\langle e_\alpha,a\rangle
\end{equation}
for a fixed normal vector $e_\alpha$.

From \eqref{2-14} for $f_i$ and  the
structure equation \eqref{2-1} , we have

\begin{equation}\label{3-3}
 f_i=\langle e_i,a\rangle.
\end{equation}
Similarly, from  \eqref{2-14} for $f_{ij}$ and the
structure equation \eqref{2-2}, we have
\begin{equation}\label{3-4}
 f_{ij}=\sum_\alpha h^\alpha_{ij}\langle e_\alpha,a\rangle.
\end{equation}

Since $a$ can be arbitrary in \eqref{3-3} and \eqref{3-4}, we
obtain (see \cite{CL})

\begin{equation}\label{3-5}
 x_i=e_i,\qquad x_{ij}=\sum_\alpha h^\alpha_{ij}e_\alpha.
 \end{equation}
Define the first derivative $g_{\alpha,i}$ of $g_\alpha$ by
\begin{equation}
\sum_ig_{\alpha,i}\theta_i=dg_\alpha+\sum_\beta
g_\beta\theta_{\beta\alpha}.
\end{equation}
We have, by use of \eqref{2-3},
\begin{equation}\label{3-7}
g_{\alpha,i}=-\sum_k h^\alpha_{ik}\langle
 e_k,a\rangle.
\end{equation}
Taking covariant derivatives on both sides of \eqref{3-7} in the
$e_j$ direction and using  \eqref{3-5}, we have

\begin{equation}\label{3-8}
 g_{\alpha,ij}=-\sum_{k}
h^\alpha_{ikj}\langle
e_k,a\rangle-\sum_{k,\beta}h^\alpha_{ik}h^\beta_{kj}\langle
e_\beta,a\rangle,
\end{equation}
where the second derivative $g_{\alpha,ij}$ of
$g_\alpha$ is defined by
\begin{equation}
\sum_jg_{\alpha,ij}\theta_j=dg_{\alpha,i}
 +\sum_{j}g_{\alpha,j}\theta_{ji}+\sum_\beta g_{\beta,i}\theta_{\beta\alpha}.
\end{equation}

Again, since  $a$ is arbitrary in \eqref{3-7} and \eqref{3-8}, we
obtain (see \cite{CL})
\begin{equation}
e_{\alpha,i}=-\sum\limits_j h^\alpha_{ij}e_j,\qquad
 e_{\alpha,ij}=-\sum\limits_{k}
h^\alpha_{ikj}e_k-\sum\limits_{k,\beta}h^\alpha_{ik}h^\beta_{kj}
e_\beta,
\end{equation}
where the covariant derivative $h^\alpha_{ijk}$ of the second
fundamental form $h^\alpha_{ij}$ is defined by \eqref{2-8}.

\vskip 3mm

 Now the self-shrinker equation \eqref{1-2} is equivalent to
\begin{equation}\label{3-11}
-H^\alpha=<x,e_\alpha>, \quad n+1\leq \alpha\leq n+p.
\end{equation}

Taking covariant derivative of \eqref{3-11} with respect to $e_i$ by
use of (3.5) and (3.10), we have
\begin{equation}\label{3-12}
-H^\alpha_{,i}=-\sum\limits_{j}h^\alpha_{ij}<x,e_j>,\quad 1\leq
i\leq n,\quad n+1\leq \alpha\leq n+p.
\end{equation}

Taking covariant derivative of (3.12) with respect to $e_k$  by use
of (3.5) and (3.11), we have
\begin{equation}
\begin{array}{lcl}
-H^\alpha_{,ik}&=&-\sum\limits_{j}h^\alpha_{ijk}<x,e_j>-h^\alpha_{ik}-\sum\limits_{\beta,j}
h^\alpha_{ij}h^\beta_{jk}<x,e_\beta>\\
&=&-\sum\limits_{j}h^\alpha_{ijk}<x,e_j>-h^\alpha_{ik}+\sum\limits_{\beta,j}
H^\beta h^\alpha_{ij}h^\beta_{jk}.
\end{array}
\end{equation}

Writing
\begin{equation}\label{3-14}
\sigma_{\alpha\beta}=\sum\limits_{i,j}h^\alpha_{ij}h^\beta_{ij},
\end{equation}
we have
\begin{equation}
\sum\limits_{\alpha,\beta}\sigma_{\alpha\beta}H^\alpha H^\beta\leq
|A|^2 |H|^2.
\end{equation}

We are now ready to prove the following key lemma:

\medskip\noindent
{\bf Lemma 3.1} {\it Let $M^n$  be an $n$-dimensional complete
self-shrinker in $\R^{n+p}$ without boundary and with polynomial
volume growth, if $|A|^2$ is bounded on $M^n$, then
\begin{eqnarray*}
\int_M|\nabla^\perp H|^2e^{-\frac{|x|^2}{2}}dv
&=&\int_M
[\sum\limits_{\alpha,\beta}\sigma_{\alpha\beta}H^\alpha H^\beta-|H|^2]e^{-\frac{|x|^2}{2}}dv\\
&\leq& \int_M [|A|^2-1]|H|^2e^{-\frac{|x|^2}{2}}dv.
\end{eqnarray*}
}

\begin{proof} Letting $i=k$ in (3.13) and summing over i, we get
\begin{equation}\label{3-16}
\Delta^\perp H^\alpha=\sum\limits_j
H^\alpha_{,j}<x,e_j>+H^\alpha-\sum\limits_\beta\sigma_{\alpha\beta}H^\beta.
\end{equation}

Since $M^n$ has polynomial volume growth and $|A|^2$ is bounded on $M^n$, \eqref{3-11}, \eqref{3-12}, \eqref{3-14} and \eqref{3-16} imply that
\begin{equation*}
\int_M|\nabla^\perp H|^2e^{-\frac{|x|^2}{2}}dv<+\infty, \qquad \int_M \sum\limits_\alpha H^\alpha \Delta^\perp  H^\alpha e^{-\frac{|x|^2}{2}}dv<+\infty,
\end{equation*}
and
\begin{equation*}
\int_M\sum_{\alpha,i}H^{\alpha}H^{\alpha}_{,i}<x,e_i>e^{-\frac{|x|^2}{2}}dv<+\infty.
\end{equation*}
Let  $\varphi_r(x)$ be a smooth cut-off function with compact support in $B_{x_0}(r+1)\subset M$,

\begin{equation*}
    \varphi_r(x)=\left\{\begin{array}{ll}
                        1, & \textrm{in }B_{x_0}(r) \\
                        0 & \textrm{in }M\setminus B_{x_0}(r+1)
                      \end{array}\right.\qquad 0\leq \varphi_r(x)\leq 1, \quad |\nabla\varphi_r|\leq 1.
\end{equation*}
Then, by integration by parts, we get
\begin{eqnarray*}
\int_M \sum\limits_\alpha \Delta^\perp H^\alpha  (\varphi_r H^\alpha) e^{-\frac{|x|^2}{2}}dv&=&\int_M \varphi_r H^{\alpha}H^{\alpha}_{,i}<x,e_i>e^{-\frac{|x|^2}{2}}dv -\int_MH^{\alpha}_{,i}(\varphi_r H^{\alpha})_{,i}e^{-\frac{|x|^2}{2}}dv\\
&=&\int_M\varphi_r \left(\sum_{\alpha,i} H^{\alpha}H^{\alpha}_{,i}<x,e_i>-|\nabla^{\perp}H|^2\right)e^{-\frac{|x|^2}{2}}dv\\
&&\quad -\int_M \sum_{\alpha,i}H^{\alpha}H^{\alpha}_{,i}(\varphi_r)_i e^{-\frac{|x|^2}{2}}dv.
\end{eqnarray*}
Letting $r\rightarrow +\infty$, the dominated convergence theorem implies that
\begin{eqnarray}\label{3-17}
\int_M \sum\limits_\alpha \Delta^\perp H^\alpha H^\alpha e^{-\frac{|x|^2}{2}}dv&=&  \int_M\left(\sum_{\alpha,i}H^{\alpha}H^{\alpha}_{,i}<x,e_i>-|\nabla^{\perp}H|^2\right)e^{-\frac{|x|^2}{2}}dv.
\end{eqnarray}
Putting \eqref{3-16} into \eqref{3-17}, we obtain:
\begin{eqnarray*}
\int_M|\nabla^\perp H|^2e^{-\frac{|x|^2}{2}}dv
&=&\int_M
\left(\sum\limits_{\alpha,\beta}\sigma_{\alpha\beta}H^\alpha H^\beta-|H|^2\right)e^{-\frac{|x|^2}{2}}dv\\
&\leq& \int_M \left(|A|^2-1\right)|H|^2e^{-\frac{|x|^2}{2}}dv.
\end{eqnarray*} \end{proof}

\smallskip\noindent
 {\bf Remark 3.1} From the proof of Lemma 3.1, one can see that the conclusion of Lemma 3.1 is valid even if $|A|^2$
 has certain growth in $|x|^2$.

\section { Proof of Theorem 1.1}

Now we present the proof of Theorem 1.1.
\begin{proof}[Proof of Theorem 1.1]

Under the assumptions of Theorem 1.1, from Lemma 3.1, we know that
either ${\bf H}\equiv 0$, or ${\bf H}\not=0$ but with $\nabla^\perp{\bf
H}\equiv 0$ and $|A|^2\equiv 1$.

If  ${\bf H}\equiv 0$, we have $<x,e_\alpha>\equiv 0,$ $n+1\leq
\alpha\leq n+p$, from which we easily conclude from (3.12) that $M$ is totally
geodesic, that is, a hyperplane in $\R^{n+1}$.

Next, suppose  that ${\bf H}\not=0$,  $\nabla^\perp{\bf H}\equiv 0$, and $|A|^2\equiv
1$.  In this case, (3.13) becomes
\begin{equation}
\sum\limits_{\beta,j} H^\beta h^\alpha_{ij}h^\beta_{jk}=h^\alpha_{ik}+
\sum\limits_{j}h^\alpha_{ijk}<x,e_j>,\quad 1\leq i,k\leq n;
n+1\leq\alpha\leq n+p.
\end{equation}

Multiplying both sides of (4.1) by $h^\alpha_{ik}$ and summing over $\alpha,i, k$, we get
\begin{equation}
\sum\limits_{\alpha,\beta,i,j,k}H^\beta
h^\alpha_{ij}h^\beta_{jk}h^\alpha_{ik}=|A|^2+\frac{1}{2}(|A|^2)_{,j}<x,e_j>
=|A|^2=1.
\end{equation}

Next we choose a local orthonormal frame  $\{e_{\alpha}\}$  for the normal
bundle of $x: M \rightarrow \R^{n+p}$, such that $e_{n+p}$ is
parallel to the mean curvature vector $\bf H$; i.e.,

\begin{equation}
e_{n+p}=\frac{\bf H}{|\bf {H}|},\quad H^{n+p}=H, \qquad H^{\alpha}=0, \quad \alpha\not=n+p.
\end{equation}

Because now the equality holds in (3.15), we have
\begin{equation}
h^{\alpha}_{ij}=0, \quad \alpha\not=n+p,\qquad
|A|^2=\sum\limits_{i,j}h^{n+p}_{ij}h^{n+p}_{ij}=1.
\end{equation}

Since $\nabla^\perp \bf{H}\equiv 0$ and $|A|^2\equiv 1$, by the
definition of  $\Delta$ and using (2.7), (2.10), (2.4),
(2.5) and (2.9), we have  (c.f. \cite{Si},\cite{SSY},\cite{Li},\cite{WA})
\begin{eqnarray*}
0&=&\f{1}{2}\Delta |A|^2\\
&=&\sum_{\alpha,i,j,k}(h^{\alpha}_{ijk})^2+\sum_{\alpha,i,j,k}h^{\alpha}_{ij}h^{\alpha}_{ijkk}\\
&=&\sum_{\alpha,i,j,k}(h^{\alpha}_{ijk})^2+\sum_{\alpha,i,j,k,m}
h^{\alpha}_{ij}h^{\alpha}_{mk}R_{mijk}+\sum_{\alpha,i,j,m}
h^{\alpha}_{ij}h^{\alpha}_{im}R_{mj}+\sum_{\alpha,\beta,i,j,k}
h^{\alpha}_{ij}h^{\beta}_{ik}R_{\beta\alpha jk}\\
&=&\sum_{\alpha,i,j,k}(h^{\alpha}_{ijk})^2+\sum_{\alpha,\beta,i,j,m}
H^{\beta}h^{\beta}_{mj}h^{\alpha}_{ij}h^{\alpha}_{im}-\sum_{\alpha,\beta,i,j,k,m}h^{\alpha}_{ij}h^{\beta}_{ij}h^{\alpha}_{mk}h^{\beta}_{mk}+2\sum_{\alpha,\beta,i,j,k}
h^{\alpha}_{ij}h^{\beta}_{ik}R_{\beta\alpha jk}.
\end{eqnarray*}
Plugging (4.2), (4.3) and (4.4) into the above identity, we conclude that
\begin{equation}
h^{\alpha}_{ijk}=0,\quad n+1\leq \alpha\leq n+p.
\end{equation}
Because $e_{n+1}\wedge_{n+2}\wedge\cdots \wedge e_{n+p-1}$ is
parallel in the normal bundle of $M$ and $h^{\alpha}_{ij}\equiv 0,
\quad \alpha\not=n+p$, by Theorem 1 of Yau \cite {Y}, we see that
$M$ is a hypersurface in  $\R^{n+1}$. So (4.5) implies that $M$ is
an isoparametric hypersurface, thus from $|A|^2=1$ we conclude that $M$
is either a round sphere $\S^n(\sqrt{n})$, or a cylinder
$\S^m(\sqrt{m})\times \R^{n-m}$, $1\leq m\leq n-1$ in $\R^{n+1}$. This
completes the proof of Theorem 1.1.
\end{proof}

\section{Further Remarks}

In this section, we make several simple observations:

\medskip\noindent
{\bf Proposition 5.1} {\it If a submanifold $M^n\to \R^{n+p}$ is an $n$-dimensional
complete self-shrinker without boundary and with polynomial volume growth, such that
\begin{equation}\label{5-1}
|H|^2\geq n,
\end{equation}
then $|H|^2\equiv n$ and $M$ is a minimal submanifold in the sphere
$\S^{n+p-1}(\sqrt{n})$.}

\begin{proof}[Proof of Proposition 5.1]  From (3.5) and (3.11), we have
\begin{equation}\label{5-2}
\frac{1}{2}\Delta |x|^2=n+<x,\Delta x>=n+\sum\limits_\alpha H^\alpha
<x,e_\alpha>=n-|H|^2
\end{equation}
Under the polynomial volume growth assumption, (1.2) and \eqref{5-2} guarantee that
\begin{equation*}
    \int_M(\Delta |x|^2)e^{-\frac{|x|^2}{2}}dv<+\infty \qquad {\mbox and} \qquad \int_M|\nabla |x|^2|^2e^{-\frac{|x|^2}{2}}dv<+\infty.
\end{equation*}
Then, by integrating by parts and the dominated convergence theorem, it follows that (similar to the proof of Lemma 3.1)
\begin{equation}\label{5-3}
\frac{1}{4}\int_M|\nabla |x|^2|^2e^{-\frac{|x|^2}{2}}dv
=\frac{1}{2}\int_M(\Delta |x|^2)e^{-\frac{|x|^2}{2}}dv
=\int_M(n-|H|^2)e^{-\frac{|x|^2}{2}}dv.
\end{equation}
From \eqref{5-1} and \eqref{5-3}, we get
$|H|^2=n$ and  $<x,x>=r^2$. Thus by (1.2) we conclude that $r=\sqrt{n}$ and $M$ is a minimal
submanifold in the sphere $\S^{n+p-1}(\sqrt{n})$.
\end{proof}

\medskip\noindent
{\bf Proposition 5.2} {\it If a submanifold $M\to \R^{n+p}$ is an $n$-dimensional compact
self-shrinker without boundary and satisfies either $|H|^2=constant$, or
\begin{equation}\label{5-4}
|H|^2\leq n,
\end{equation}
then $|H|^2\equiv n$ and $M$ is a minimal submanifold in the sphere
$\S^{n+p-1}(\sqrt{n})$.}

\begin{proof}[Proof of Proposition 5.2] Integrating \eqref{5-2} over $M$ and using the Stokes theorem,
we have
\begin{equation}\label{5-5}
\int_M(n-|H|^2)dv=0.
\end{equation}
Hence Proposition 5.2 follows from \eqref{5-5}, \eqref{5-4}, and (1.2).
\end{proof}

\medskip\noindent
{\bf Remark 5.1} Let $x:M\to \R^{n+p}$ be an $n$-dimensional
submanifold. If $x$ satisfies
\begin{equation}
\lambda H^\alpha=<x,e_\alpha>, \quad n+1\leq \alpha\leq n+p
\end{equation}
for some positive constant $\lambda$, then we call $M$ a {\it
self-expander} of the mean curvature flow. Observe that for a self-expander, we have
\begin{equation}\label{5-7}
\frac{1}{2}\Delta |x|^2=n+<x,\Delta x>=n+n\sum\limits_\alpha
H^\alpha <x,e_\alpha>=n+n\lambda |H|^2.
\end{equation}
From \eqref{5-7}, we immediately get

\medskip\noindent
{\bf Proposition 5.3} {\it There exists no n-dimensional compact
self-expander without boundary in $\R^{n+p}$.}


\medskip

Finally, we list some simple examples of self-shrinkers of higher codimensions.

\medskip\noindent
{\bf Example 5.1} {\it  For any positive integers $m_1,\cdots,m_p$ such that $m_1+\cdots
+m_p=n$, the submanifold
\begin{equation}
M^n=\S^{m_1}(\sqrt{m_1})\times\cdots\times\S^{m_p}(\sqrt{m_p})\subset\R^{n+p}
\end{equation}
is an n-dimensional compact self-shrinker in $\R^{n+p}$ with
\begin{equation}
{\bf H}=-X,\qquad |{\bf H}|^2=n,\qquad |A|^2=p
\end{equation} Here
\begin{equation}
\S^{m_i}(r_i)=\{X_i\in\R^{m_i+1}: |X_i|^2=r_i^2\},\qquad
i=1,\cdots,p
\end{equation}
is a $m_i$-dimensional round sphere with radius $r_i$.}

\medskip\noindent
{\bf Example 5.2} {\it For positive integers $m_1,\cdots,m_p,q\geq 1$, with
$m_1+\cdots+m_p+q=n$, the submanifold
\begin{equation}
M^n=\S^{m_1}(\sqrt{m_1})\times\cdots\times\S^{m_p}(\sqrt{m_p})\times\R^q\subset\R^{n+p}
\end{equation}
is an n-dimensional complete non-compact self-shrinker in $\R^{n+p}$
with polynomial volume growth which satisfies
\begin{equation}
{\bf H}=-X^{\perp},\qquad |{\bf H}|^2=\sum_{i=1}^pm_i,\qquad |A|^2=p.
\end{equation}
}

\medskip\noindent
{\bf Remark 5.2} In Example 5.1 and Example 5.2, if we let $p\geq 2$, then we have
an n-dimensional self-shrinker of codimension $p$ with $|A|^2=p\geq 2$, thus not
one of the three cases in Theorem 1.1.

\smallskip\noindent
 {\bf Remark 5.3} From Example 5.2, we can see that the condition ``$|{\bf H}|^2\geq n$'' in Proposition
5.1 is necessary.

\medskip\noindent
{\bf Example 5.3} (cf. \cite{CA}) {\it Let
\begin{equation}
X: \S^2(\sqrt{m(m+1)})\hookrightarrow
\S^{2m}(\sqrt{2})\subset\R^{2m+1},\qquad m\geq 2
\end{equation}
be a minimal surface in $\S^{2m}(\sqrt{2})$. Consider it as a surface in $\R^{2m+1}$, then it
is a self-shrinker with
\begin{equation}
{\bf H}=-X,\qquad |{\bf H}|^2=2,\qquad |A|^2=2-\f {2}{m(m+1)}<2,
\end{equation}}

\smallskip\noindent
 {\bf Remark 5.4} By choosing local orthogonal frame  $\{e_{\alpha}\}$  for the normal
bundle of $x: M^n \rightarrow \R^{n+p}$, such that $e_{n+p}$ is
parallel to the mean curvature vector $\bf H$, by Lemma 3.1, if
$|A|^2$ is bounded, and
 \begin{equation}
 \sum_{i,j}h_{ij}^{n+p}h_{ij}^{n+p}\leq 1,\label{5.15}
 \end{equation}
 we have $\nabla^{\bot}{\bf H}=0$, that is,  $|{\bf H}|^2=constant$ and unit mean curvature vector
 field $\nu={\bf H}/|{\bf H}|$ is parallel
in the normal bundle. From Proposition 5.2 and Theorem 1.3 of Smoczyk \cite{Smo},
 we have

\medskip\noindent
{\bf Proposition 5.4} {\it  Let $M^n$  be an $n$-dimensional
complete self-shrinker in $\R^{n+p}$ without boundary and with
polynomial volume growth. If $|A|^2$ is bounded on $M^n$ and
(5.15) holds, then
$$
M^n=N^m\times \R^{n-m}, \qquad 0\leq m\leq n, \label{5.15}
$$
where  $N^m$ is a $m$-dimensional minimal submanifold in $\S^{m+p-1}(\sqrt{m})$.}

\begin{flushleft}
\medskip\noindent
\begin{tabbing}
XXXXXXXXXXXXXXXXXXXXXXXXXX*\=\kill
Huai-Dong Cao\>Haizhong Li\\
Department of Mathematics\>Department of Mathematical Sciences\\
Lehigh University\>Tsinghua University\\
Bethlehem, PA 18015\>100084, Beijing\\
USA\>People's Republic of China\\
E-mail:huc2@lehigh.edu\>E-mail:hli@math.tsinghua.edu.cn
\end{tabbing}

\end{flushleft}

\end{document}